   \def\B{\cal B}
   
\def\W{{\cal W}}

\def\Zet{{\Bbb Z}}

\def\k{{\bold k}}
\def\B{\cal B}

\def\Bn{B_{n}}
\def\Pn{P_{n}}
\def\Vk{\b{V}}

\def\TV{TV}
\def\Va{V^{\ast}}
\def\Vd{V^{\otimes 2}}
\def\Vtr{V^{\otimes 3}}

\def\sg{\sigma}
\def\tr{\triangle}

\def\bT{bTV}

\def\Hom{\operatorname{Hom}}

\def\End{\operatorname{End}}

\def\group{\operatorname{group}}
\def\alg{\operatorname{alg}}
\def\coalg{\operatorname{coalg}}

\def\Vt{V^{\otimes}}

\newcommand{\id}{\operatorname{id}}

\newtheorem{lem}{Lemma}[section]
\newtheorem{th}[lem]{Theorem}

\documentstyle[12pt]{amsart}

\begin{document}

\title{Braided antisymmetrizer\\ as bialgebra homomorphism}
\author{Jerzy R\'o\.za\'nski} \thanks {Supported by KBN grant
\# 2 P302 023 07}
\address{Department of Mathematical Methods in Physics,
University of Warsaw, Ho\.za 74, PL 00-682 Warszawa, Poland}
\email{rozanski@@fuw.edu.pl}
\date{December 1995}
\maketitle

\begin{abstract}
For an Yang Baxter operator we show that a bialgebra
homomorphism from a free braided tensor bialgebra
to a cofree braided shuffle bialgebra
is the Woronowicz braided antisymmetrizer \cite{SW}.
A cofree braided shuffle bialgebra is a braided generalization of a
cofree shuffle bialgebra introduced by Sweedler \cite{Sweedler}.
Its graded dual bialgebra is a free braided tensor bialgebra \cite{OPR}.
\end{abstract}

\section{Introduction}
Let $\k$ be an associative ring with unit and $V$ be a $\k$-bimodule.
For an Yang Baxter operator $\B\in\End (V^{\otimes 2})$ Wo\-ro\-no\-wicz
\cite{SW} defined a braided antisymmetrizer $W(\B)\in\End (\Vt)$,
which is the gene\-ra\-lization of the antisymmetrizer
for the flip $\B=-\tau$.

For the operator $W_{n}(\B)=W(\B)|V^{\otimes n}$ exist two operators
$\mu (\B)$ and $\tr (\B)$ such that
$$
W_{k+l}(\B)=\mu_{k,\;l}(\B)\circ \{W_{k}(\B)\otimes W_{l}(\B)\}
=\{W_{k}(\B)\otimes W_{l}(\B)\}\circ\tr_{k,\;l} (\B).
$$

We point out that two braided bialgebras are associated
with those decompositions,
the cofree braided shuffle bialgebra $bShV$ and its graded dual: the free
braided tensor bialgebra $\bT$. We show that the cofree braided shuffle
bialgebra is the generalization of the cofree shuffle bialgebra considered by
Sweedler \cite{Sweedler}.

We prove that the kernel of the braided antisymmetrizer $W(\B)$ is the
biideal in the free braided bialgebra $\bT$.

\noindent{\bf Theorem}.\\
{\em The homomorphism $W$ from the
braided tensor bialgebra $\bT$ to the braided shuffle bialgebra $bShV$,
such that $W|\k=\id$ and $W|V=\id$, is the braided antisymmetrizer $W(\B)$.}

\section{Notations}
Let $\Bn$ be a braid group with generators
$\{\sg_{1},\sg_{2},\ldots,\sg_{n-1}\}$ and relations
\begin{gather}\label{brarel}
\sg_{i}\sg_{j}=\sg_{j}\sg_{i}\qquad |i-j| > 1,\nonumber\\
\sg_{i}\sg_{i+1}\sg_{i}=\sg_{i+1}\sg_{i}\sg_{i+1}\qquad i=1,\;2,\ldots,\;n-2.
\end{gather}

The generators of a permutation group $\Pn$: $\{t_{1},t_{2},\ldots,t_{n-1}\}$
satisfy relations (\ref{brarel}) and the additional relation
$\forall i=1,\ldots,n-1:\;t_{i}^{2}=1$.
For a permutation $p\in\Pn$ let $I(p)$ be the number of inversed pairs
in the sequence $(p(1),\;p(2),\ldots,\;p(n))$, then
\begin{equation}\label{perm}
p=t_{k_1}t_{k_2}\cdots t_{k_{I(p)}}.
\end{equation}

One can define \cite{SW} a map $\pi$ from $\Pn$ to a braid group $\Bn$.
For the permutation (\ref{perm}) we define $\pi (t_{i})=\sg_{i}$ and
\begin{equation}\label{elemb}
b=\pi (p)=\sg_{k_1}\sg_{k_2}\cdots\sg_{I(p)}\in\Bn.
\end{equation}

A braid (\ref{elemb}) is independent of the choice of the decomposition
(\ref{perm}), see \cite{SW}.
Denote by $\Xi$ the image of the map $\pi$,
\begin{equation}
\Xi_{n}:=\pi (\Pn)\subset \Bn.
\end{equation}

A $(k,\;l)$-shuffle for $k,\;l\geq 0$ is a permutation $p\in P_{k+l}$
satisfying
\begin{equation}
p(1)<p(2)<\ldots<p(k)\;\text{and}\;
p(k+1)<p(k+2)<\ldots<p(l).
\end{equation}
A $(0,\;n)$ and $(n,\;0)$ shuffle is an identity.
The subset of $(k,\;l)$-shuffles will be denoted by
$Sh_{k,\;l}\subset P_{k+l}$. For this set we have the corresponding
subset of braided shuffles
\begin{equation}
bSh_{k,\;l}:=\pi (Sh_{k,\;l})\subset \Bn.
\end{equation}

For $\k$-bimodule $V$ denote by $\Vk$ the $\Zet$-graded $\k$-bimodule
\begin{equation}
\Vk=\k\oplus V\oplus \Vd\oplus\cdots .
\end{equation}
Let $[v_{1}| v_{2}|,\ldots,|v_{n}]$ be an element of
$\Vk^{n}\equiv\Vk\cap V^{\otimes n}$ and $v_{0}=1\in \k$.
The tensor algebra $\TV$ means in this case \cite{GS} the pair
$\TV=\{\Vk,\;\otimes\}$, where
$\otimes: \Vk \times \Vk \rightarrow \Vk$ is a free product.

Let $\Va=\Hom (V,\; \k)$ be the dual $\k$-bimodule for $V$.
For the evaluation\\
$<\Va,\; V>$ are two possibilities to evaluate $<\Va\otimes\Va,\;V\otimes V>.$
We use the following evaluation
(note that transposition depends on this choice),
 \begin{equation}
<\alpha\otimes\beta,\;v\otimes w>=<\beta,\; v><\alpha,\; w>,\qquad
\forall \alpha,\;\beta\in\Va\quad \text{and}\quad \forall v,\;w\in V.
\end{equation}

For the bimodule $\Vk$ the graded dual bimodule $\Vk^{g}$ means
\cite{Sweedler}
\begin{equation}
\Vk^{g}=\k\oplus V^{\ast}\oplus {V^{\ast}}^{\otimes 2}\oplus \ldots .
\end{equation}

Let $\B$ be the Yang Baxter operator, i.e. an invertible endomorphism of the
$\k$-bimodule $\Vd$ which satisfies the braid equation
\begin{equation}
(\id\otimes \B)(\B\otimes \id)(\id\otimes \B)=
(\B\otimes\id)(\id\otimes \B)(\B\otimes\id)\in\End(\Vtr).
\end{equation}
Let the set
$B_{k}=\id_{k-1}\otimes B\otimes\id_{n-k-1}$ of
$V^{n}\subset\Vk\in\End (\Vk^{\otimes n})$ be bimodule endomorphisms.
For the Yang Baxter operator $\B$ we can define the representation
$\rho_{\B}$ of the braid group $B_{n}$ in the $\Zet$-graded
$\k$-bimodule $\Vk$,
\begin{equation}\label{rrep}
\rho_{\B}\in\group\{B_{n},\;\End(\Vk^{n})\}:\qquad
\rho_{\B}(\sg_{k})=B_{k}.
\end{equation}

\section{Two braided bialgebras}
For the Yang Baxter operator $\B$ on the bimodule $V$ the subset
$\Xi_{n}\subset\Bn$ is acting on the $\Zet$-graded $\k$-bimodule
through the representation $\rho_{\B}$ (\ref{rrep}).
Let us introduce the notion of the cofree braided shuffle algebra.
For the $\k$-bimodule $\Vk$ one can define \cite{Fox} the cofree
comultiplication $\tr^{\otimes}$
\begin{equation}
\tr^{\otimes} [v_{1}| v_{2}|\ldots|v_{n}] =\sum_{k=0}^{n}
[v_{1}| v_{2}|\ldots |v_{k}] \otimes [v_{k+1}|\ldots |v_{n}].
\end{equation}

The braided shuffle multiplication $\mu (\B):\Vk\otimes\Vk\rightarrow\Vk$
is defined for $\mu_{k,n-k}(\B)=\mu (\B)|\Vk^{k}\otimes\Vk^{n-k}$
\begin{gather}\label{mu}
\mu_{k,n-k} (\B)=\sum_{b\in bSh_{k,n-k}} \rho_{\B}(b):
\Vk^{k}\otimes \Vk^{n-k}\rightarrow \Vk^{n},\nonumber\\
\mu_{k,n-k} (\B)([v_{1}| v_{2}|\ldots |v_{k}]
\otimes [v_{k+1}|\ldots |v_{n}])
=\sum_{b\in bSh_{k,n-k}} \rho_{\B} [v_{1}| v_{2}|\ldots |v_{n}].
\end{gather}

For example $\mu_{1,2} (\B)=\id_{3}+\B\otimes\id+
(\B\otimes\id)(\id\otimes\B)$.
\begin{lem}
The multiplication $\mu (\B)$ (\ref{mu}) is associative and $\B$-braided.
\end{lem}
\begin{pf}
Denote an element $[v_{1}| v_{2}|\ldots |v_{i}]$ by $v_{I}$. Then we have
$$
\mu (\B) (v_{I}\otimes v_{J})=\sum_{K=I+J} v_{K}\in \Vk^{K}.
$$
The associativity condition is proved by the following equation
\begin{multline*}
[\mu (\B)\circ(\mu (\B)\otimes \id)] (v_{I}\otimes v_{J}\otimes v_{K})\\
=\sum_{L=I+J+K} v_{L}
=[\mu (\B)\circ(\id\otimes \mu (\B))] (v_{I}\otimes v_{J}\otimes v_{K}).
\end{multline*}\end{pf}

Consider the $\B$-braided monoidal category. A bialgebra is defined
over this category if the multiplication $m$ and
comultiplication $\tr$ are morphisms and satisfy the following
compatibility condition
\begin{gather}
\tr\circ m=(m\otimes m)\circ (\id\otimes\B\otimes\id)\circ (\tr\otimes\tr ).
\end{gather}

\begin{lem}\label{bb}
The triple $(\Vk,\;\mu (\B),\;\tr^{\otimes})$ is the braided bialgebra.
\end{lem}
The proof by induction for the term $\Vk^{k}\otimes \Vk^{l}$ is omitted.

Consider the $\Zet$-graded dual $\k$-bimodule $V^{g}$.
The free braided tensor bialgebra $\bT$ is the graded dual to the cofree
braided shuffle bialgebra. Then the free braided tensor bialgebra $\bT$
means the triple $\{\Vk^{g},\;\otimes,\;\tr^{\mu} (\B)\}$, where the
multiplication $\otimes$ and the comultiplication $\tr^{\mu}$ are graded
duals in the following sense
\begin{gather}
\tr^{\mu} (\B)={\mu}^{g} (\B),\qquad \text{and}
\qquad \otimes={\tr^{\otimes}}^{g}.
\end{gather}
Dual version of the lemma \ref{bb} is the following braided bialgebra.
\begin{lem}
The triple $(\Vk,\;\otimes,\;\tr^{\mu} (\B))$ is the braided bialgebra.
\end{lem}

\section{Braided antisymmetrizer}
For an Yang Baxter operator $\B$ the braided antisymmetrizer $W(\B)$ was
defined by Woronowicz \cite{SW} as
\begin{equation}
W(\B)=\sum_{b\in\Xi_{n}} \rho_{\B}(b).
\end{equation}

For $\B=-\pi$ we get the Woronowicz form of the braided antisymmetrizer
with the sign of the permutation.

\begin{lem}
For the braided antisymmetrizer $\W (\B)$,
the multiplication $\mu (\B)$ and the comultiplication $\tr^{\mu} (\B)$
the following recurrent formula is holds
$$W_{n+1}(\B)=\mu(\B)_{n,\;1}\circ [W_{n}(\B)\otimes\id]
=[\id\otimes W_{n}(\B)]\circ {\tr^{\mu}}_{1,\;n} (\B).$$
\end{lem}
For the proof see \cite{SW} for details.

{}From this lemma by induction we can prove
\begin{equation}\label{Wdec}
W_{k+l}(\B)=\mu(\B)_{k,\;l}\circ [W_{k}(\B)\otimes W_{l}(\B)].
\end{equation}
and dually
\begin{equation}
W_{k+l}(\B)=[W_{k}(\B)\otimes W_{l}(\B)]\circ {\tr^{\mu}_{1,\;n}(\B)}.
\end{equation}

Consider braided bialgebras $\bT$ and $bShV$.
A map $W: \bT\rightarrow BShV$ is a homomorphism of these bialgebras
$W\in bialg(\bT,\; bShV)$ if satisfies two conditions
\begin{itemize}
\item W is the algebra map:
$$W\in\alg (\otimes,\;\mu (\B)),\qquad
W\circ\otimes=\mu (\B)\circ (W\otimes W)
$$
\item W is the coalgebra map:
$$W\in\coalg (\tr^{\mu} (\B),\;\tr^{\otimes}),\qquad
\tr^{\otimes}\circ W=(W\otimes W)\circ\tr^{\mu} (\B).$$
\end{itemize}

\begin{th}\label{main}
The bialgebra homomorphism $W\in bialg(\bT,\; bShV)$, such that $W|\k=\id$
and $W|V=\id$, is the braided antisymmetrizer W(\B).
\end{th}
\begin{pf}
{}From the assumption that W is the algebra map we obtain
$$W_{n}=\mu_{n-1,1} (\B)\circ (W_{n-1}\otimes \id).$$
Then by induction we have
$$
W_{n}=\mu_{n-1,1} (\B)\circ\mu_{n-2,1} (\B)\circ\ldots\circ \mu_{1,1} (\B).
$$
{}From the assumption $W|\k=\id$ and $W|V=\id$ this is the braided
antisymmetrizer $$W=W(\B).$$
\end{pf}
\section*{Acknowledgments}
I would like to thank prof. Z. Oziewicz for fruitful discussions.

\end{document}